\documentclass[12pt,reqno]{amsart}

\usepackage{amsmath,amsthm,amsfonts,amscd,flafter,epsf}

\let\qua\stdspace

\newtheorem{theorem}{Theorem}[section]

\newtheorem{cor}[theorem]{Corollary}
\newtheorem{prop}[theorem]{Proposition}

\newtheorem{remark}[theorem]{Remark}

\newcommand{\HF}{HF}

\newcommand{\Q}{\mathbb{Q}}

\newcommand{\Z}{\mathbb{Z}}

\newcommand{\Image}{\mathrm{Im}}

\newcommand{\cm}{\cdot}

\newcommand{\ModSWfour}{\mathcal{M}}
\newcommand{\ModFlow}{\ModSWfour}

\newcommand{\SpinC}{{\mathrm{Spin}}^c}
\newcommand{\Spin}{{\mathrm{Spin}}}

\newcommand\abuts\Rightarrow
\newcommand\Sym{\mathrm{Sym}}

\newcommand\Eul{\widehat{\chi}}

\newcommand\HFpRed{\HFp_{\red}}

\newcommand\RelSpinC{\underline{\SpinC}}

\newcommand\ModSphere{\ModFlow\left({\mathbb S}\longrightarrow
\Sym^{g-1}(\Sigma_{1})\times \Sym^2(\Sigma_{2})\right)}
\newcommand\ModSpheres\ModSphere

\newcommand\HFpred{\HFp_{\rm red}}

\newcommand\HFred{\HF_{\rm red}}

\newcommand{\red}{\mathrm{red}}

\newcommand\HFp{\HFb}

\newcommand\HFinf{HF^\infty}

\newcommand\HFb{HF^+}

\newcommand\UnparModSp{\widehat \ModSp}
\newcommand\UnparModFlow\UnparModSp
\newcommand\Mod\ModSp

\newcommand{\spinc}{\mathfrak s}

\newcommand{\spinct}{\mathfrak t}

\newcommand\ModMaps{\mathcal M}
\newcommand\ModSp\ModMaps

\newcommand\Dual{\mathcal D}
\newcommand\Duality\Dual

\begin{document}
\title[The Renormalized Euler Characteristic and $L$-space Surgeries]{The Renormalized Euler Characteristic \\and $L$-space Surgeries}
\author{Raif Rustamov}
\address{The Program in Applied and Computational Mathematics, Princeton University\\New Jersey 08540, USA}
\email{rustamov@princeton.edu}
\begin{abstract}
Using the equivalence between the renormalized Euler
characteristic of Ozsv\'ath and Szab\'o, and the Turaev torsion
normalized by the Casson-Walker invariant, we make calculations
for $S^3_{p/q}(K)$. An alternative proof of a theorem by Ozsv\'ath
and Szab\'o on $L$-space surgery obstructions is provided.
\end{abstract}
\maketitle

\section{Introduction}

In the beautiful paper \cite{unk}, Ozsv\'ath and Szab\'o consider
the \emph{correction terms} of $Y=S^3_{2n-1/2}(K)$, where $K$ is a
knot in $S^3$, and prove that if $Y$ is an $L$-space, then there
is a symmetry among these correction terms. Given an $L$-space
whose correction terms are known, one can check whether this
symmetry is satisfied, and if not, it would follow that the
manifold at hand cannot be obtained as a $\frac{2n-1}{2}$ surgery
on a knot in $S^3$. This obstruction is then used to calculate
some new unknotting numbers. The symmetry was later generalized to
the case of $p/q$ surgery in \cite{RationalSurgeries}.

We remind that for a rational homology sphere $Y$ and a $\SpinC$
structure $\spinct$ on it, the \emph{renormalized Euler
characteristic} $\Eul(Y,\spinct)$ is defined as
$$\Eul(Y, \spinct) =
\chi(\HFred(Y,\spinct)) - \frac{1}{2} d(Y,\spinct),$$ where
$d(Y,\spinct)$ denotes the correction terms. By definition, $Y$ is
an $L$-spaces if and only if $\HFpred(Y,\spinct) \cong 0$. One
immediately sees that for $L$-spaces any symmetry of correction
terms is equivalent to the symmetry of $\Eul$ 's. This points out
that calculating the renormalized Euler characteristic for
$S^3_{p/q}(K)$ could be interesting.

Given a knot $K$, let
$$\Delta_K(T)=a_0+\sum_{j>0}
a_j\left(T^j+T^{-j}\right)$$ be the symmetrized Alexander
polynomial of $K$, and define $$t_i(K) = \sum_{j \geq 1}
ja_{|i|+j}.$$
\begin{theorem}
Let $p\neq 0$ and $q$ be relatively prime integers. For each
$\spinct$, the renormalized Euler characteristic
$\Eul(S^3_{p/q}(K), \spinct)$ can be expressed in terms of the
coefficients of the Alexander polynomial of $K$.
\end{theorem}
This theorem  follows from the identification
$\Eul=-\tau+\lambda$, where $\tau$ is the Turaev torsion, and
$\lambda$ is the Casson-Walker invariant, see \cite{rustamov}. In
fact, we give a precise formula for the renormalized Euler
characteristic in Proposition \ref{slprop}. The next theorem is
proved using this formula.
\begin{theorem}
Let $p$ and $q$ be relatively prime integers with $p/q>1$. Suppose
that the Alexander polynomial of knot $K$ satisfies $a_j=0$ for
$j>\frac{p}{2q}+1$. Let $n$ be any integer, and denote $r=\lceil
\frac{pn}{q}-1\rceil \in \Z/p\Z$, then for any $|i|\leq p/2q$ we
have
$$\Eul(S^3_{p/q}(K),r+i)-\Eul(S^3_{p/q}(U),r+i)=t_i,$$
where $U$ is the unknot, and we have used the affine
identification $\SpinC(S^3_{p/q}(K)) \cong \Z/p\Z$ explained in
Section~\ref{prelim}; also see Remark \ref{spins}.
\end{theorem}

Note that if $p/q\leq 1$, then the condition of the theorem will
force $a_j=0$ for $j>1$. In this case, it follows from the
calculations that
$$\Eul(S^3_{p/q}(K),-1)-\Eul(S^3_{p/q}(U),-1)=\lceil q/p\rceil a_1,$$
and the difference is zero for all remaining $\SpinC$ structures.\

As a corollary of this theorem we obtain an alternative proof for
the Theorem 1.2 of \cite{RationalSurgeries}.
\begin{cor} Let $K$ be a knot which admits an $L$-space surgery
for some $\frac{p}{q}>1$. Then, for all integers $i$ with $|i|
\leq \frac{p}{2q}$ we have that
$$d(S^3_{p/q}(K),i-1)-d(S^3_{p/q}(U),i-1)=-2t_i,$$
while for all $|i| > \frac{p}{2q}$, we have $t_i(K) = 0$.
\end{cor}
\begin{proof}
If $S^3_{p/q}(K)$ is an $L$-space, then it is known that
$S^3_{\lceil p/q\rceil}(K)$ is also an $L$-space. However, then
$\lceil p/q\rceil \geq 2g(K)-1$, where $g(K)$ is the degree of the
Alexander polynomial of $K$. Thus, $a_j=0$ for $j>\frac{\lceil
p/q\rceil+1}{2}$, i.e. also for $j>\frac{p}{2q}+1$. The corollary
follows by taking $n=0$ in the previous theorem and observing that
for $L$-spaces
$$\Eul(Y, \spinct) = - \frac{1}{2} d(Y,\spinct).$$
\end{proof}

The organization of this paper is as follows: the required
preliminaries are presented in Section 2. Calculations needed for
Theorem 1.1  are made in Section 3. We finish with the proof of
Theorem 1.2 in Section 4.

\medskip
\noindent{\bf{Acknowledgments}}\qua I am pleased to thank my
advisor Zolt\'an Szab\'o for very helpful conversations.

\section{Preliminaries}
\label{prelim} Let $Y$ be a rational homology sphere, $\spinct$ be
a $\SpinC$ structure on it. We can consider the Heegaard Floer
homology group $\HFp(Y,\spinct)$. This is a $\Q$ graded module
over $Z[U]$. We can also consider a simpler version,
$\HFinf(Y,\spinct)$ for which one can prove
\begin{equation}
\label{eq:HFinf} \HFinf(Y,\spinct) \cong Z[U,U^{-1}]
\end{equation}
 for each $\spinct$. There is a natural $Z[U]$ equivariant map  $$\pi \colon\HFinf(Y,\spinct)\longrightarrow \HFp(Y,\spinct)$$ which is zero in sufficiently negative degrees and an isomorphism in all sufficiently positive degrees.
$\HFpred(Y, \spinct)$ is defined as
$$\HFpRed(Y, \spinct)=\HFp(Y,\spinct)/ \Image \pi.$$
Let $d(Y,\spinct)$ be the \emph{correction term}, defined as the
minimal degree of any non-torsion class of $\HFp(Y, \spinct)$
lying in the image of $\pi$. As was mentioned, the
\emph{renormalized Euler characteristic} $\Eul(Y,\spinct)$ is
defined by
$$\Eul(Y, \spinct) = \chi(\HFred(Y,\spinct)) - \frac{1}{2} d(Y,\spinct).$$

The \emph{Casson-Walker} invariant, $$\lambda : \{\mathrm{rational
\qua homology \qua spheres \qua modulo \qua homeomorphisms}\}
\rightarrow \Q,$$ is an extension of Casson's invariant to
rational homology spheres, see \cite{Walker}. It perhaps worth
noting that in our normalization $\lambda(\Sigma(2,3,5))=-1$,
where $\Sigma(2,3,5)$ is oriented as the boundary of the negative
definite $E_8$ plumbing.

We will also consider the normalized Casson-Walker invariant given
by $\lambda'(Y) = \big|H_1(Y;\Z)\big|\lambda(Y)$. Note that
$\lambda'(S^3) = 0$. Let $U$ be the unknot, and $K$ be any knot in
$S^3$. Walker's surgery formula implies that
\begin{eqnarray}
\label{walker}\lambda'(S^3_{p/q}(K)) - \lambda'(S^3_{p/q}(U)) =q
\cm \sum_{j\geq 1} j^2 a_j ,
\end{eqnarray}
where $a_j$ are the coefficients of the Alexander polynomial of
$K$.

The \emph {Turaev torsion function} is another important invariant
of three-manifolds. Turaev first defined it for combinatorial
Euler structures, but later the connection with these and $\SpinC$
structures emerged, providing us with a function
$$\tau(Y,\cdot):\SpinC(Y) \rightarrow \Z,$$
see \cite{TuraevSurgery}.

The Turaev torsion can be also introduced in the case of
three-manifolds with torus boundary. Indeed, consider the manifold
$M=S^3-\mathrm{nd}(K)$ with torus boundary. We denote the set of
\emph{relative} $\SpinC$ structures on M by $\RelSpinC(M)$. The
Turaev torsion of this manifold is a function
$$\tau(M,\cdot):\RelSpinC(M) \rightarrow \Z.$$
Actually, under a certain affine map $\RelSpinC(M)\cong
\Z_{\mathrm{odd}}$, one has
$$\tau(M, k)=\mathrm{sign}(k)\cm\sum_{j\geq\frac{|k|+1}{2}} a_j,$$
for any odd integer $k$, see 9.1.4 of \cite{TuraevSurgery}. Here
$a_i$ are again the coefficients of the Alexander polynomial of
$K$. In what follows, we will use this particular identification
of $\RelSpinC(M)$ with $\Z_{\mathrm{odd}}$. Of course, given a
relative $\SpinC$ structure on $M=S^3-\mathrm{nd}(K)$ we can look
at the $\SpinC$ structure on $S^3_{p/q}(K)$ which extends it.
Thus, we can get an affine identification $\SpinC(S^3_{p/q}(K))
\cong \Z/p\Z$  as follows: if $k \in \Z_{\mathrm{odd}} \cong
\RelSpinC(M)$ extends to $i \in \Z/p\Z \cong
\SpinC(S^3_{p/q}(K))$, then
$$i\equiv\frac{k-1}{2} \mod p.$$

Let $x$ be any integer that satisfies $qx \equiv -1 \mod p$.
Remember that $H_1(S^3_{p/q}(K);\Z)$ is generated by the meridian
of the knot $K$, and the homology class of $K$ is $x$ times this
generator. The following formula is a consequence of 10.6.3.2 of
\cite{Turaev}:
\begin{align}
\label{tursurg}\nonumber\tau(S^3_{p/q}(K), i) &-
\tau(S^3_{p/q}(K), i+x)=\\&= \tau(S^3_{p/q}(U), i) -
\tau(S^3_{p/q}(U), i+x)+\sum_{\{k | i\equiv\frac{k-1}{2} \mod p\}}
\tau(M, k),
\end{align}
where $i \in \Z/p\Z$. Note that, to get this formula we have
replaced the purely homological data of the original formula by
terms coming from the surgery on unknot.

Let $Y$ be a rational homology sphere, $\spinct$ be $\SpinC$
structure on it.  We have the following formula for the
renormalized Euler characteristic, see \cite{rustamov}:
\begin{equation}
\label{myformula}
 \Eul(Y, \spinct) =-\tau(Y,\spinct) +\lambda(Y).
\end{equation}
Note that since Turaev torsions add up to zero, we have
\begin{equation}
\label{sum} \sum_{\spinct\in\SpinC(Y)}\Eul(Y,\spinct)=\lambda'(M).
\end{equation}
\section{Calculations}
From now on we fix the knot $K$, and two relatively prime integers
$p$ and $q$. Note that calculating the renormalized Euler
characteristics $\Eul(S^3_{p/q}(K),i)$ is equivalent to computing
$$S_i=\Eul(S^3_{p/q}(K),i)-\Eul(S^3_{p/q}(U),i),$$
because the second term is already known, \cite{AbsGraded}. For $0
\leq i <p$ let
$$T_i=\sum_{\{k | i\equiv\frac{k-1}{2} \mod p\}} \tau(M, k),$$
which is readily calculated from the Alexander polynomial of $K$.
Combining (\ref{tursurg}) and (\ref{myformula}) we get
$$S_{i+x}-S_{i}=T_i,$$
for any $i \in \Z/p\Z$. Fix any $l \in \Z/p\Z$. To get a formula
for $S_l$, we write the following equations
$$S_{l+(j+1)x}-S_{l+jx}=T_{l+jx},$$
for $j=0,1,...,p-1$. It follows that $S_{l+x}=S_{l}+T_l$,
$S_{l+2x}=S_l+T_l+T_{l+x}$  and so on. Using (\ref{walker}) and
(\ref{sum}), we get $$p S_l + \sum_{j=0}^{p-1}(p-j-1)T_{l+jx} = q
\cm \sum_{j\geq 1} j^2 a_j.$$ Thus,
\begin{equation}
\label{sl} S_l =  \frac{1}{p} \left( q\cm \sum_{j\geq 1} j^2
a_j-\sum_{j=0}^{p-1}(p-j-1)T_{l+jx}\right),
\end{equation}
 which establishes the
following proposition.
\begin{prop}
\label{slprop} Let $K$ be a knot in $S^3$, $p>0$ and $q$
relatively prime integers, and $x=-q^{-1}\mod p$, then
$$\Eul(S^3_{p/q}(K),l)=\Eul(S^3_{p/q}(U),l)+\frac{q}{p} \cm
\sum_{j\geq 1} j^2
a_j-\frac{1}{p}\cm\sum_{j=0}^{p-1}(p-j-1)T_{l+jx},$$ where
$\Eul(S^3_{p/q}(U),i) =\Eul(L(-p,q),i) = -d(L(-p,q),i)/2.$
\end{prop}
\begin{remark}
\label{spins} Let $p$ be odd, so there is a unique $\Spin$
structure on $S^3_{p/q}(K)$. In the canonical affine
identification $\SpinC(S^3_{p/q}(K)) \cong \Z/p\Z$ this $\Spin$
structure corresponds to $0$, and conjugation is equivalent to
multiplication by $-1$. We denote the $\SpinC$ structure
corresponding to $i \in \Z/p\Z$ in this identification by
$\spinc_i$ (thus, $\spinc_0$ is the $\Spin$ structure, $\bar
\spinc_i=\spinc_{-i}$), and in our identification by $\spinct_i$.
Let us spell out the correspondence between these two
identifications. Note that we always have $T_{\frac{p-1}{2}} =0$,
which means $S_{\frac{p-1}{2}+x}=S_{\frac{p-1}{2}}$. This
universal equality must be a consequence of the conjugation
symmetry $\Eul(Y, \bar\spinct)=\Eul(Y,\spinct)$. Thus,
$\spinct_{\frac{p-1}{2}+x}$ and $\spinct_{\frac{p-1}{2}}$ are
conjugate. If $\spinc_i =\spinct_{\frac{p-1}{2}}$, then $i+x=-i
\mod p$, i.e. $i=\frac{p-1}{2}x \mod p$. As a result,
$\spinct_{\frac{(p-1)(1-x)}{2} \mod p}$ is the $\Spin$ structure
of $S^3_{p/q}(K)$.
\end{remark}

\section{Proof of Theorem 1.2}
From now on we assume that $q>1$. The case when $q=1$ is only
slightly different, and is left to the reader. Note that we can
simplify the previous formula under certain conditions on the
Alexander polynomial. Let the knot $K$ be such that $a_j=0$ for
$j\geq p/2$, which means that for $0\leq i<p/2$,$$T_i=\sum_{\{k |
i\equiv\frac{k-1}{2} \mod p\}} \tau(M, k)
=\tau(M,2i+1)=\sum_{j\geq i+1} a_j,$$ and similarly, for $p/2\leq
i\leq p-1$,
$$T_i=-\sum_{j\geq p-i} a_j.$$
With $l$ fixed, let $0 \leq u_j\leq p-1$ satisfy $$l+u_jx=j-1 \mod
p,$$ and $0 \leq v_j\leq p-1$ satisfy $$l+v_jx=-j \mod p.$$ Let
$$c_i=\sum_{j=1}^{i} (u_j-v_j),$$ and since $x=-q^{-1} \mod
p$, one has
$$c_i=p\cm\sum_{j=1}^{i} \left(\left\{\frac{q(l+1-j)}{p}\right\}
-\left\{\frac{q(l+j)}{p}\right\}\right),$$ where
$\{\alpha\}=\alpha-\lfloor\alpha\rfloor$ denotes the fractional
value of the number $\alpha$. Now we can rewrite (\ref{sl}) as
\begin{equation}
\label{simplified} S_l =  \frac{1}{p}  \sum_{i\geq 1} (qi^2+c_i)
a_i.
\end{equation}

The proof of Theorem 1.2 now becomes an exercise in arithmetic.
The assumption of the theorem is that $p/q>1$ and $q>1$ , and that
the Alexander polynomial of knot $K$ satisfies $a_j=0$ for
$j>\frac{p}{2q}+1$. Given any $l \in \Z/p\Z$, let us find the
coefficient of $a_1$ in $S_l$. We have find the value of
$$c_1=p\left(\left\{\frac{ql}{p}\right\}
-\left\{\frac{q(l+1)}{p}\right\}\right).$$Obviously,
$\frac{q(l+1)}{p} - \frac{ql}{p} =\frac{q}{p}<1$, which means that
$c_1=-q$ unless there is an integer $n$ so that
$$\frac{ql}{p}<n
\leq  \frac{q(l+1)}{p},$$ i.e. $l=\lceil\frac{qn}{p} -1\rceil$, in
which case we have $c_1=p-q$. As a result, the coefficient of
$a_1$ in $S_l$ is equal to $1$ if $l=\lceil\frac{qn}{p} -1\rceil$,
and equal to $0$ otherwise.

Any $a_j$ with $j<\frac{p}{2q}+\frac{1}{2}$ can be analyzed
similarly, with the result that $a_j$'s coefficient in $S_l$ is
equal to zero unless $l=\lceil\frac{qn}{p} -1\rceil + i$, where
$|i|\leq j-1$, in which case the coefficient is equal to $j-i$.

The case of $a_j$ with $\frac{p}{2q}+\frac{1}{2}\leq j \leq
\frac{p}{2q}+1$ is a bit different, because the difference
$\frac{q(l+j)}{p} - \frac{q(l+1-j)}{p}$ is not necessarily less
than $1$, thus there is a chance that between these two fractions
two integers may appear. However, since we are interested in $S_l$
with $l=\lceil\frac{qn}{p} -1\rceil + i$ , where $i\leq p/2q$,
this does not happen. As a result, $a_j$'s coefficient in $S_l$
for $l=\lceil\frac{qn}{p} -1\rceil + i$, where $|i|\leq j-1$ is
equal to $j-i$. Since $a_j=0$ for $j>\frac{p}{2q}+1$, we have
proved that given any integer $n$, for every $|i|\leq p/2q$ we
have
$$\Eul(S^3_{p/q}(K),r+i)-\Eul(S^3_{p/q}(U),r+i)=t_i,$$
where $r=\lceil \frac{pn}{q}-1\rceil \in \Z/p\Z$.

\end{document}